\documentstyle[amssymb,10pt]{amsart}
\newtheorem{theorem}{Theorem}[section]
\newtheorem{lemma}[theorem]{Lemma}
\newtheorem{proposition}[theorem]{Proposition}
\newtheorem{corollary}[theorem]{Corollary} 
\theoremstyle{definition}  
\newtheorem{definition}[theorem]{Definition}

\newtheorem{remark}[theorem]{Remark}

\newcommand{\id}{\operatorname{id}}

\newcommand{\Ad}{\operatorname{Ad}}

\newcommand{\nc}{\newcommand}
\nc{\Symm}{{\on{Sym}}}

\newcommand{\on}{\operatorname}   
\newcommand{\eps}{\varepsilon}
 \nc{\cE}{{\cal E}}

\renewcommand{\O}{{\mathcal{O}}}

\renewcommand{\a}{{\mathfrak a}}

\renewcommand{\c}{{\mathfrak c}}

\nc{\SL}{{\mathfrak sl}}
\nc{\HH}{{\mathfrak h}}
\newcommand{\g}{{\mathfrak{g}}}

\newcommand{\m}{{\mathfrak{m}}}

\nc{\wh}{\widehat}\nc{\wt}{\widetilde}

\newcommand{\ben}{\begin{enumerate}}
\newcommand{\een}{\end{enumerate}}

\newcommand{\cO}{{\mathcal O}}

\newcommand{\cR}{{\mathcal R}}

\newcommand{\RR}{{\bold R}}
\newcommand{\LL}{{\bold L}}

\newcommand{\ZZ}{{\mathbb{Z}}}
\newcommand{\KM}{{\mathbb{K}}}

\renewcommand{\Ad}{{\operatorname{Ad}}}

\begin{document}

\title[Coboundary Lie bialgebras and commuting subalgebras]
{Coboundary Lie bialgebras and commutative subalgebras of 
universal enveloping algebras}

\begin{abstract} 
We solve a functional version of the problem of twist 
quantization of a coboundary Lie bialgebra $(\g,r,Z)$. 
We derive from this the following results: 
(a) the formal Poisson manifolds $\g^*$ and $G^*$ are 
isomorphic; (b) we construct an injective algebra morphism 
$S(\g^*)^\g \hookrightarrow U(\g^*)$. When $(\g,r,Z)$ 
can be quantized, we construct a deformation of this morphism. 
In the particular case when $\g$ is quasitriangular and nondegenerate, 
we compare our construction with Semenov-Tian-Shansky's construction 
of a commutative subalgebra of $U(\g^*)$. We also show that the canonical 
derivation of the function ring of $G^*$ is Hamiltonian.  
\end{abstract}

\author{Benjamin Enriquez}
\address{IRMA (CNRS), rue Ren\'e Descartes, F-67084 Strasbourg, France}
\email{enriquez@@math.u-strasbg.fr}

\author{Gilles Halbout}
\address{IRMA (CNRS), rue Ren\'e Descartes, F-67084 Strasbourg, France}
\email{halbout@@math.u-strasbg.fr}

\maketitle

\centerline{0. {\textsc{Introduction}}}

\bigskip

Let $(\g,r,Z)$ be a coboundary Lie bialgebra over a field $\KM$
of characteristic $0$. This means that $\g$ is a Lie bialgebra, 
the Lie cobracket $\delta$ of which is the coboundary of 
$r\in\wedge^2(\g)$: 
$\delta(x)=[x \otimes 1 + 1 \otimes x,r]$ for any $x \in \g$. 
This condition means that $Z:= \on{CYB}(r)$ belongs to 
$\wedge^3(\g)^\g$
(here $\on{CYB}$ is the l.h.s. of the classical Yang-Baxter
equation). Quasi-triangular and triangular Lie bialgebras 
are particular cases of this definition. 

\medskip

It is an open question to construct a twist quantization of $(\g,r,Z)$, 
i.e., a pair $(J,\Phi)$, where $J \in U(\g)^{\otimes 2}[[\hbar]]$
and $\Phi\in U(\g)^{\otimes 3}[[\hbar]]$ are invertible ($\hbar$ is a formal
series), $\Phi$ is $\g$-invariant, $(J,\Phi)$ satisfies a cocycle relation
and deforms $(r,Z)$.
If $(J,\Phi)$ satisfies these conditions, then $\Phi$ satisfies the 
pentagon relation, is $\g$-invariant and deforms $Z$: 
such a $\Phi$ is called a 
quantization of $Z$. In \cite{Dr:QH}, Proposition 3.10, Drinfeld constructed 
a quantization $\Phi$ of $Z$. Any pair $(J,\Phi)$ can be made admissible 
(in the sense of \cite{EH2}) and the associated formal functions 
$(\rho,\varphi)$ 
then satisfy functional analogues of the pentagon and cocycle equations
(this is explained in Section \ref{sect:rel}).
We call this system of equations the functional analogue of twist quantization.   

\medskip 

We describe the set of solutions of the functional analogue of twist 
quantization for $(\g,r,Z)$. Namely, we derive from Drinfeld's result 
that $Z$ can be lifted to an element $\varphi \in \m_{\g^*}^{\wh\otimes 3}$ 
satisfying the functional pentagon relation (Proposition \ref{lift1}). 
We then prove that $r$ can be lifted to an element 
$\rho \in \m_{\g^*}^{\wh\otimes 2}$, such that $(\rho,\varphi)$ 
satisfies the functional cocycle relation (Theorem \ref{lift2}); 
here $\m_{\g^*} \subset \cO_{\g^*}$ is the maximal ideal of 
the ring of formal functions on $\g^*$. We show that all solutions are 
related by suitable gauge transformations. 

\smallskip

The first corollary is that the formal Poisson manifolds $\g^*$ and 
$G^*$ are isomorphic (Corollary \ref{coroprinc}).

\smallskip

In Section 3, we prove another corollary (Theorem \ref{theoprinc}): 
we construct an injection of algebras $S(\g^*)^\g \hookrightarrow 
U(\g^*)$ (here $\g^*$ acts on $S(\g)$ by symmetric powers of the 
coadjoint action). This morphism is filtered, its associated 
graded is the canonical inclusion
$S(\g^*)^\g \subset S(\g^*)$. This way, we obtain a commutative subalgebra of 
$U(\g^*)$. The fact that the graded subalgebra $S(\g^*)^\g \subset S(\g^*)$ is 
Poisson commutative can be viewed as a classical limit of this situation.
It can either be viewed as a corollary of the fact that  
$\cO_G^\g \subset \cO_G$ is a Poisson commutative subalgebra
or it can be proved directly (Lemma \ref{lemma:poisson}); 
here $\cO_G$ is the ring of formal functions on $G$, on which $\g$ acts by
conjugation.   

\smallskip

In Sections \ref{sect:O} and \ref{sect:SG}, assuming the existence 
of a twist quantization of $\g$, we 
construct formal deformations of the algebra inclusions 
$S(\g^*)^\g \hookrightarrow U(\g^*)$ and $\cO_G^\g\subset \cO_{G^*}$. 
All these results use the theory of duality of QUE and QFSH algebras; 
this theory is recalled in Section \ref{sect:duality}.  

\smallskip

In Section \ref{sect:Sem}, we assume that $\g$ is quasitriangular. 
In that case, we show that $U(\g^*)$ contains a family 
$C_s$ of commutative subalgebras, indexed by $s\in \KM$; 
this result may be viewed as a classical limit of Drinfeld's 
result about commutativity of twisted traces. We explain why 
only $C_0$ has an analogue in the general coboundary case. In 
\cite{STS1}, Semenov-Tian-Shansky defined an algebra morphism 
$U(\g)^\g \stackrel{\on{STS}}{\to} U(\g^*)$; we show that its image 
coincides with $C_1$, and is therefore in general different from 
the image $C_0$ of the morphism in our construction. 

\smallskip 
Finally, in Section \ref{sect:D}, we show that the canonical 
derivation of $\cO_{G^*}$ is Hamiltonian. This derivation is equal to  
$\hbar^{-1}(S^2-\id)_{|\hbar=0}$, where $S$ is the antipode of any 
quantization of $\cO_{G^*}$. 

\subsection*{Notation}

We use the standard notation for the coproduct-insertion maps:
we say that an ordered set is a pair of a finite set $S$ and a 
bijection $\{1,\ldots,|S|\} \to S$. 
For $I_1,\dots,I_m$ disjoint ordered subsets of $\{1,\dots,n\}$,
$(U,\Delta)$ a Hopf algebra and $a \in  U^{\otimes m}$,
we define
$$a^{I_1,\dots,I_n}= \sigma_{I_1,\ldots,I_m} \circ 
(\Delta^{|I_1|}\otimes \cdots \otimes \Delta^{|I_n|})(a),
$$
with $\Delta^{(1)}=\on{id}$, $\Delta^{(2)}=\Delta$, 
$\Delta^{(n+1)}=({\on{id}}^{\otimes n-1} \otimes \Delta)\circ \Delta^{(n)}$, 
and $\sigma_{I_1,\ldots,I_m} : U^{\otimes \sum_i |I_i|} \to 
U^{\otimes n}$ is the 
morphism corresponding to the map $\{1,\ldots,\sum_i |I_i|\} 
\to \{1,\ldots,n\}$ taking $(1,\ldots,|I_1|)$ to $I_1$, 
$(|I_1| + 1,\ldots,|I_1| + |I_2|)$ to $I_2$, etc. 
When $U$ is cocommutative, this definition depends only on
the sets underlying $I_1,\ldots,I_m$.

\subsection*{Acknowledgements}
We would like to thank V. Dolgushev, P. Etingof and L.-C. Li for
discussions. 

\section{Solutions of the functional twist equations}

If $\g$ is a Lie algebra, we denote by 
$\cO_{\g^*} = \wh S(\g)$ the formal series ring of functions on the 
formal neighborhood of $0$ in $\g^*$. We define by $\m_{\g^*} \subset 
\cO_{\g^*}$ the maximal ideal of this ring. If $k$ is an integer $\geq 1$, 
we denote by $\cO_{(\g^*)^k} = \wh S(\g)^{\wh\otimes k}$ the
\footnote{$\wh\otimes$ is 
the completed tensor product, defined by 
$V_0[[x_1,\ldots,x_n]] \wh\otimes W_0[[y_1,\ldots,y_n]] := 
V_0\otimes W_0[[x_1,\ldots,y_n]]$, where $V_0,W_0$ are vector spaces.} 
ring of formal functions functions on $(\g^*)^k$, by $\m_{(\g^*)^k}$
its maximal ideal and by $\m_{(\g^*)^k}^i$ the $i$th power of this ideal. 

If $f,g\in \m_{(\g^*)^k}^2$, then the series 
$f \star g = f + g + {1\over 2} \{f,g\} + 
\cdots + B_n(f,g) + \cdots$ is convergent, where $\sum_{i\geq 1} 
B_i(x,y)$ is the Baker-Campbell-Hausdorff series specialized to the 
Poisson bracket of $\m_{(\g^*)^k}^2$. The product $\star$ defines a 
group structure on $\m_{(\g^*)^k}^2$. 

If $f\in \cO_{\g^*}^{\wh\otimes n}$ 
and $P_1,\dots, P_m$ are disjoint subsets of 
$\{1,\dots,m\}$, one defines $f^{P_1,\dots,P_n}$ as
in the Introduction using the cocommutative coproduct of 
$\cO_{\g^*}$ (dual to the addition of $\g^*$).

\medskip

Let $\g$ be a Lie algebra and $Z\in \wedge^3(\g)^\g$. 

\begin{proposition}
\label{lift1}
There exists $\varphi\in (\m_{\g^*}^{\wh\otimes 3})^\g 
(\subset \m_{(\g^*)^3}^2)$, the image of 
which under the map 
$\m_{\g^*}^{\wh\otimes 3} \to 
(\m_{\g^*} / \m_{\g^*}^2)^{\otimes 3} = \g^{\otimes 3} 
\stackrel{\on{Alt}}{\to} \wedge^3(\g)$
equals $Z$ (here $\on{Alt}$ is the total antisymmetrization map)
and satisfying the functional pentagon equation  
$$
\varphi^{1,2,34} \star \varphi^{12,3,4} = \varphi^{2,3,4} \star 
\varphi^{1,23,4} \star \varphi^{1,2,3}. 
$$
Such a $\varphi$ (we call it a lift of $Z$) is unique up to the action 
of an element of $(\m_{\g^*}^{\wh\otimes 2})^\g$ 
by $\sigma \cdot \varphi = \sigma^{2,3} \star \sigma^{1,23} \star 
\varphi \star (-\sigma)^{12,3} \star (-\sigma)^{1,2}$. 
\end{proposition}

{\em Proof.} In \cite{Dr:QH}, Proposition 3.10, Drinfeld constructed a solution 
$\Phi\in U(\g)^{\otimes 3}[[\hbar]]$ of the pentagon equation 
\begin{equation} \label{pent}
\Phi^{1,2,34} \Phi^{12,3,4} = \Phi^{2,3,4} \Phi^{1,23,4} \Phi^{1,2,3} 
\end{equation}
such that $\eps^{(2)}(\Phi)=1$ and $\Phi = 1^{\otimes 3} + O(\hbar)$
(here $\eps^{(2)} = \id\otimes \eps\otimes\id$; applying $\eps$ to the 
first and third factors of (\ref{pent}), we also get $\eps^{(1)}(\Phi) =
\eps^{(3)}(\Phi)=1$). 

In \cite{EH2}, we stated that $\Phi$ can be transformed into an admissible 
solution $\Phi'$ of the same equations, using an invariant twist. In Appendix 
\ref{app:A}, we explain why the proof given in \cite{EH2} is wrong and
we give a correct proof.  

The classical limit of $\hbar\log(\Phi')$ then satisfies the functional 
pentagon equation. 
This gives the existence of $\varphi$. One can also construct $\varphi$ 
directly using cohomological methods, as it will be done for $\rho$ later. 

Let us prove uniqueness: let $\varphi$ and $\varphi'$ be two lifts of
$Z$. The classes of $\varphi$ and $\varphi'$ are the same in 
$\m_{\g^*}^{\wh\otimes 3}/(\m_{\g^*}^{\wh\otimes 3}\cap \m_{(\g^*)^3} ^3)$, 
as this space is $0$.
Let $N$ be an integer $\geq 3$; assume that we have found $\sigma_N
\in (\m_{\g^*}^{\wh\otimes 2})^\g$ such that $\sigma_N \cdot \varphi$ 
and $\varphi'$ are equal modulo $\m_{\g^*}^{\wh\otimes 3} \cap 
\m_{(\g^*)^3}^N$. Write $\varphi'=\sigma_N \cdot \varphi+\psi$,
with $\psi \in (\m_{\g^*}^{\wh\otimes 3} \cap \m_{(\g^*)^3}^N)^\g$.
We will use the following lemma (see \cite{EGH}, p. 2477): 
\begin{lemma}
\label{lemmetech}
For any $k\geq 1$ and $n \geq 2$, 
$f,h \in \m_{(\g^*)^k}^2$ and 
$g \in \m_{(\g^*)^k}^n$, one has
$$
f \star (h + g) = f \star h + g, \quad 
(f+g)\star h = f \star h +g
\hbox{\ modulo\ }\m_{(\g^*)^k}^{n+1}.
$$ 
\end{lemma}

Let $\overline\psi$ be the class of $\psi$ in 
$(\m_{\g^*}^{\wh\otimes 3}\cap\m_{(\g^*)^3}^N)^\g /
(\m_{\g^*}^{\wh\otimes 3}\cap \m_{(\g^*)^3}^{N+1})^\g = 
(S^{>0}(\g)^{\otimes 3})^\g_N$. Then
$\overline\psi^{1,2,34} +\overline\psi^{12,3,4} = 
\overline\psi^{2,3,4} + \overline\psi^{1,23,4} + \overline\psi^{1,2,3}$,
which means that $\overline\psi$ is a cocycle in the subcomplex 
$((S^{>0}(\g)^{\otimes \cdot})^\g,d)$ of the 
co-Hochschild\footnote{We 
denote by $S(\g)$ the symmetric algebra of $\g$, by $S^{>0}(\g)$
is positive degree part; the index $N$ means the part of total degree $N$.} 
complex $(S(\g)^{\otimes \cdot},d)$. 
Using \cite{Dr:QH}, Proposition 3.11, one can prove that the $k$th
cohomology group of this complex is $\wedge^k(\g)^\g$ and that the
antisymmetrization map coincides with the canonical map from the space of 
cocycles to the cohomology.
For $N=3$, the hypothesis implies that $\on{Alt}(\overline\psi) = 0$, 
so $\overline\psi$ is a coboundary of an element $\overline\tau_3 \in 
(S^{>0}(\g)^{\otimes 2})^\g_3$. For $N >3$, $\overline\psi$ is the 
a coboundary of an element $\overline\tau_N\in (S^{>0}(\g)^{\otimes 2})^\g_N$, 
since the degree $N$ part of the relevant cohomology group 
vanishes. 
We then set $\sigma_{N+1} = \sigma_N + \tau_N$, where $\tau_N
\in (\m_{\g^*}^{\wh\otimes 2} \cap \m_{(\g^*)^2}^N)^\g$ is a lift of 
$\overline\tau_N$. Then 
$\sigma_{N+1} \cdot \varphi$ and $\varphi'$ are equal modulo 
$\m_{\g^*}^{\wh\otimes 3} \cap \m_{(\g^*)^3}^{N+1}$. The sequence 
$(\sigma_N)_{N\geq 3}$ has a limit $\sigma$. Then $\sigma \cdot \varphi =
\varphi'$. 
\hfill \qed \medskip 

We now construct a lift of $r$: 
\begin{theorem}
\label{lift2}
There exists $\rho\in \m_{\g^*}^{\wh\otimes 2}$, the image of which in 
$\g^{\otimes 2}$ under the square of 
the projection $\m_{\g^*} \to \m_{\g^*} / \m_{\g^*}^2 = \g$
equals $r$, and such that 
\begin{equation}
\label{twistequation}
\rho^{1,2}\star \rho^{12,3} = \rho^{2,3} \star \rho^{1,23} \star \varphi. 
\end{equation}
Such a $\rho$ (we call it a lift of $r$) is unique up to the action of 
$\m_{\g^*}$ 
by $\lambda \cdot \rho = \lambda^{1} \star \lambda^{2} \star \rho 
\star (-\lambda)^{12}$. We call equation (\ref{twistequation}) the functional 
cocycle equation. 
\end{theorem}

{\em Proof.}
Let us construct $\rho$ by induction: we will construct a convergent sequence 
$\rho_N \in \m_{\g^*}^{\wh\otimes 2}$ ($N\geq 2$)
satisfying (\ref{twistequation}) in $\m_{\g^*}^{\wh\otimes 3} /
(\m_{\g^*}^{\wh\otimes 3} \cap \m_{(\g^*)^3}^N)$. When $N = 3$, we take for  
$\rho_2$ any lift of $r$ to $\m_{\g^*}^{\wh\otimes 2}$; then 
equation (\ref{twistequation}) is automatically satisfied. 

Let $N$ be an integer $\geq 3$; 
assume that we have constructed $\rho_N$ in $\m_{\g^*}^{\wh\otimes 2}$
satisfying equation (\ref{twistequation}) in $\m_{\g^*}^{\wh\otimes 3}/
(\m_{\g^*}^{\wh\otimes 3}\cap\m_{(\g^*)^3} ^N)$.
Set $\alpha_N := 
\rho_N^{1,2}\star \rho_N^{12,3} -\rho_N^{2,3} \star \rho_N^{1,23} 
\star \varphi$. Then $\alpha_N$ belongs to  
$\m_{\g^*}^{\wh\otimes 3}\cap\m_{(\g^*)^3}^N$, and the following 
equalities hold in $\m_{\g^*}^{\wh\otimes 4}/(\m_{\g^*}^{\wh\otimes 4} \cap 
\m_{(\g^*)^4}^{N+1})$:  
\begin{align*}
\alpha_N^{12,3,4} = &~ 
\rho_N^{1,2} \star \alpha_N^{12,3,4}= 
\rho_N^{1,2} \star \rho_N^{12,3}\star \rho_N^{123,4} - 
\rho_N^{1,2} \star \rho_N^{3,4} \star \rho_N^{12,34} \star \varphi^{12,3,4} 
\\ & ~(\hbox{using Lemma \ref{lemmetech}})
\\ 
=&~(\alpha_N^{1,2,3}+ \rho_N^{2,3} \star \rho_N^{1,23}\star \varphi^{1,2,3})
\star \rho_N^{123,4}-\rho_N^{3,4}\star \rho_N^{1,2}\star \rho_N^{12,34} 
\star \varphi^{12,3,4} \\
=&~\alpha_N^{1,2,3} + \rho_N^{2,3} \star \rho_N^{1,23}
\star \rho_N^{123,4}\star \varphi^{1,2,3}
-\rho_N^{3,4}\star (\rho_N^{2,34}\star \rho_N^{1,234} \star \varphi^{1,2,34}
+\alpha_N^{1,2,34}) \star \varphi^{12,3,4} \cr
&~(\hbox{using Lemma \ref{lemmetech}, the invariance of }\varphi\hbox{ and
the definition of }\alpha_N^{1,2,34})\cr
=&~\alpha_N^{1,2,3} + \rho_N^{2,3} \star (\alpha_N^{1,23,4}+
\rho_N^{23,4} \star \rho_N^{1,234} \star \varphi^{1,23,4})\star 
\varphi^{1,2,3}\cr
&- \alpha_N^{1,2,34}-\rho_N^{3,4}\star \rho_N^{2,34}\star \rho_N^{1,234} 
\star \varphi^{1,2,34} \star \varphi^{12,3,4} \cr
&~(\hbox{using the definition of }\alpha_N^{1,23,4}
\hbox{ and Lemma \ref{lemmetech}})\cr
=&~\alpha_N^{1,2,3} + \alpha_N^{1,23,4}+ (\rho_N^{3,4} \star 
\rho_N^{2,34}\star \varphi^{2,3,4}+ \alpha_N^{2,3,4} )\star \rho_N^{1,234} 
\star \varphi^{1,23,4}\star \varphi^{1,2,3}\cr
&- \alpha_N^{1,2,34}-\rho_N^{3,4}\star \rho_N^{2,34}\star \rho_N^{1,234} 
\star \varphi^{1,2,34} \star \varphi^{12,3,4} \cr
&~ (\hbox{using the definition of }\alpha_N^{2,3,4} \hbox{ and Lemma 
\ref{lemmetech}}) \cr
&~ = \alpha_N^{1,2,3} + \alpha_N^{1,23,4} - \alpha_N^{1,2,34} + \alpha_N^{2,3,4} 
\end{align*}
(using Lemma \ref{lemmetech}, the invariance of $\varphi$ and the 
fact that $\varphi$ satisfies the functional pentagon equation). 

Let us denote by $\overline\alpha_N$ the image of $\alpha_N$ in 
$(\m_{\g^*}^{\wh\otimes 3}\cap\m_{(\g^*)^3}^N) / 
(\m_{\g^*}^{\wh\otimes 3}\cap\m_{(\g^*)^3}^{N+1})
= (S^{>0}(\g)^{\otimes 3})_N$, then we get 
$$
\overline\alpha_N^{12,3,4} + \overline\alpha_N^{1,2,34}
=\overline\alpha_N^{1,2,3} + \overline\alpha_N^{1,23,4}
+\overline\alpha_N^{2,3,4}.
$$
This means that $\alpha$ is a cocycle for the subcomplex 
$(S^{>0}(\g)^{\otimes\cdot},d)$ of the co-Hochschild complex.
Using \cite{Dr:QH}, Proposition 3.11, one proves that the $k$th 
cohomology group of this subcomplex is $\wedge^k(\g)$, and that
antisymmetrization coincides with the canonical projection from the space of
cocycles to the cohomology group.
For $N=3$, the equation $\on{CYB}(r)=Z$ implies 
$\on{Alt}(\overline\alpha_3)=0$, hence $\overline\alpha_3$ is the coboundary 
of an element $\overline\beta_3\in (S^{>0}(\g)^{\otimes 2})_3$. 
For $N>3$, $\overline\alpha_N$ is the coboundary of an element 
$\overline\beta_N \in (S^{>0}(\g)^{\otimes 2})_N$, since the degree $N$ 
part of the cohomology vanishes. We then set $\rho_{N+1} := \rho_N + \beta_N$, 
where $\beta_N\in\m_{\g^*}^{\wh\otimes 2} \cap \m_{(\g^*)^2}^N$ is a
representative of $\overline\beta_N$. Then $\rho_{N+1}$ satisfies 
(\ref{twistequation}) in $\m_{\g^*}^{\wh\otimes 3}
/(\m_{\g^*}^{\wh\otimes 3} \cap \m_{(\g^*)^3}^{N+1})$. 

The sequence $(\rho_N)_{N\geq 2}$ has a limit $\rho$, which then satisfies
(\ref{twistequation}). 

The second part of the theorem can be proved either by analyzing the 
choices for $\overline\beta_N$ in the above proof, or following the proof 
of the previous proposition.
\hfill \qed \medskip 

\begin{remark} \label{rem:sigma} 
If $\varphi$ is replaced by $\varphi' = \sigma \star \varphi$, then a
solution of (\ref{twistequation}) is $\rho' = \rho \star (-\sigma)$.  
\end{remark}

\section{Isomorphism of formal Poisson manifolds $\g^* \simeq G^*$}

Let us assume that $\g$ is a finite dimensional coboundary Lie bialgebra. 
the following result was proved in \cite{EEM} when $\g$ is quasitriangular; 
the result of \cite{EEM} is itself a generalization of the formal version 
of the Ginzburg-Weinstein isomorphism (\cite{GW,A,Bo}).

\begin{corollary}
\label{coroprinc}
There exists an isomorphism of formal Poisson manifolds $\g^* \simeq G^*$.
\end{corollary}

{\em Proof.} Let $P : \wedge^2(\cO_{\g^*}) \to \cO_{\g^*}$ be the Poisson
bracket on $\cO_{\g^*}$ corresponding to the 
Lie-Poisson\footnote{or Kostant-Kirillov-Souriau, or linear} Poisson 
structure on $\g^*$. 
Then $(\O_{\g^*},m_0,P,\Delta_0)$ is a Poisson formal series Hopf (PFSH) 
algebra; it corresponds to the formal Poisson-Lie group $(\g^*,+)$
equipped with its Lie-Poisson structure. 

Set ${}^{\rho}\Delta(f) = \rho\star\Delta_0(f)\star (-\rho)$ for any 
$f\in\cO_{\g^*}$. It follows from the fact that $\rho$ satisfies the
functional cocycle equation that $(\cO_{\g^*},m_0,P,{}^{\rho}\Delta_0)$
is a PFSH algebra. 

Let us denote by ${\bf PFSHA}$ and ${\bf LBA}$ the categories of PSFH algebras
and Lie bialgebras. We have a category equivalence 
$c : {\bf PFSHA} \to {\bf LBA}$, taking $(\cO,m,P,\Delta)$ to the Lie bialgebra
$(\c,\mu,\delta)$, where $\c := \m/\m^2$ ($\m\subset\cO$ is the maximal ideal), 
the Lie cobracket of $\c$ is induced by $\Delta - \Delta^{2,1} : \m\to 
\wedge^2(\m)$, and the Lie bracket of $\c$ is induced by the Poisson 
bracket $P : \wedge^2(\m) \to \m$. The inverse of the functor $c$
takes $(\c,\mu,\delta)$ to $\cO = \wh S(\c)$ equipped with its usual product; 
$\Delta$ depends only on $\delta$ and $P$ depends on $(\mu,\delta)$. 

Then $c$ restricts to a category equivalence $c_{\on{fd}} : 
{\bf PFSHA}_{\on{fd}} \to {\bf LBA}_{\on{fd}}$ of subcategories of 
finite-dimensional objects (in the case of ${\bf PFSH}$, we say that $\cO$ is
finite-dimensional iff $\m/\m^2$ is). 

Let $\on{dual} : {\bf LBA}_{\on{fd}} \to {\bf LBA}_{\on{fd}}$ be the duality
functor. It is a category antiequivalence; we have $\on{dual}(\g,\mu,\delta) =
(\g^*,\delta^t,\mu^t)$. Then $\on{dual} \circ c_{\on{fd}} : 
{\bf PFSHA}_{\on{fd}} \to {\bf LBA}_{\on{fd}}$ is a category antiequivalence.
Its inverse it the usual functor $\g\mapsto U(\g)^*$. If $G$ is the formal 
Poisson-Lie group with Lie bialgebra $\g$, one sets $\cO_G = U(\g)^*$. 

Let us apply the functor $c$ to $(\cO_{\g^*},m_0,P,{}^\rho\Delta_0)$. 
We obtain $\c = \m/\m^2 = \g$; the Lie bracket is unchanged w.r.t. 
the case $\rho=0$, so it is the Lie bracket of $\g$; the Lie cobracket 
is given by $\delta(x) = [r,x\otimes 1 + 1\otimes x]$ since the reduction of 
$\rho$ modulo $(\m_{\g^*})^2\wh\otimes \m_{\g^*} + \m_{\g^*} \wh\otimes 
(\m_{\g^*})^2$ is equal to $r$. 

Then applying $\on{dual} \circ c_{\on{fd}}$ to 
$(\cO_{\g^*},m_0,P,{}^\rho\Delta_0)$, we obtain the Lie bialgebra 
$\g^*$. So this PFSH algebra is isomorphic to the PFSH algebra of the formal 
Poisson-Lie group $G^*$. In particular, the Poisson algebras 
$\cO_{\g^*}$ and $\cO_{G^*}$ are isomorphic. It is easy to check that 
the map $\g = \m_{\g^*} / \m_{\g^*}^2 \to \m_{G^*}/\m_{G^*}^2 = \g$
induced by this isomorphism is the identity (here $\m_{G^*} \subset \cO_{G^*}$
is the maximal ideal).  
\hfill \qed \medskip 

\begin{remark} When $\g$ is infinite dimensional, one can define 
$\cO_{G^*}$ as the image of $\g$ under ${\bf LBA} \to {\bf PFSHA}$
and then show that the Poisson algebras $\cO_{G^*}$ and $\cO_{\g^*} = 
(\wh S(\g),$ linear Poisson structure)
are isomorphic. 
\end{remark}

\section{The morphism $S(\g^*)^\g \hookrightarrow U(\g^*)$}

In this section, $\g$ is a finite dimensional coboundary Lie bialgebras.  
The following fact is well-known (\cite{STS2}): 

\begin{lemma} \label{lemma:O}
$\cO_G^\g \subset \cO_G$ is a Poisson commutative subalgebra. 
\end{lemma}

Here the action of $\g$ on $\cO_G$ corresponds to adjoint
action of $G$. We recall the proof: if $f,g\in \cO_G$, then 
$\{f,g\} = m(({\bf L} - {\bf R})(r)(f\otimes g))$, where 
${\bf L}, {\bf R}$ are the infinitesimal left and right actions
and $m$ is the product map. 
If $\varphi \in \cO_G^\g$, then ${\bf L}(a)(\varphi) = {\bf R}(a)(\varphi)$
for any $a\in\g$, 
therefore if $f,g\in\cO_G^\g$, then $({\bf L} - {\bf R})(r)(f\otimes g) = 0$, 
hence $\{f,g\}=0$.  

The inclusion $\cO_G^\g \subset \cO_G$ is a morphism of Poisson 
algebras with a decreasing filtration. By passing to the associated 
graded, we obtain: 

\begin{lemma} \label{lemma:poisson}
$S(\g^*)^\g \subset S(\g^*)$ is a Poisson commutative subalgebra. 
\end{lemma}

{\em Another proof.} If $\alpha,\beta\in\g^*$, then $[\alpha,\beta] = 
\on{ad}^*(R(\beta))(\alpha) - \on{ad}^*(R(\alpha))(\beta)$, 
where $R : \g^* \mapsto \g$ is given by $R(\xi) = (\on{id}\otimes\xi)(r)$. 

Let $f,g \in S(\g^*)^\g$ be of degrees $k$ and $\ell$.  
Write $f = \sum_\alpha a_1^\alpha \cdots a_k^\alpha$, 
$g = \sum_\beta b_1^\beta \cdots b_\ell^\beta$. Then
$$ 
\{f,g\} =  
\sum_\beta \sum_{j=1}^\ell
b_1^\beta \cdots \check b_j^\beta\cdots b_\ell^\beta
\on{ad}^*(R(b_j^\beta))(f)  
- \sum_\alpha \sum_{i=1}^k
a_1^\alpha \cdots \check a_i^\alpha\cdots a_k^\alpha
\on{ad}^*(R(a_i^\alpha))(g). 
$$  
When $f$ and $g$ are both invariant, this bracket vanishes. 
\hfill \qed \medskip 

We now prove that $S(\g^*)^\g \subset S(\g^*)$ is also the 
associated graded of an inclusion of noncommutative algebras 
with an increasing filtration:
 
\begin{theorem}
\label{theoprinc}
There exists a morphism of filtered algebras:
$$\theta :  S(\g^*)^\g \to U(\g^*),$$
the associated graded morphism of which is the canonical inclusion 
$S(\g^*)^\g \subset S(\g^*)$.
\end{theorem}

{\em Proof.} Let us denote by ${\bf FSHA}$ the category of 
formal series Hopf (FSH) algebras and by ${\bf FilAlg}$ the category 
of filtered
algebras. There is a  contravariant functor (restricted duality) 
${\bf FSHA} \to {\bf FilAlg}$, defined by $\cO\mapsto \cO^\circ$, 
where $\cO^\circ = \{\ell\in \cO^* | \exists n\geq 0, \ell(\m^n) = 0\}
\subset \cO^*$; here $\m\subset\cO$ is the maximal ideal of $\cO$.
The algebra structure of $\cO^\circ$ is defined by $(\ell_1 \cdot \ell_2)(f) 
= (\ell_1 \otimes \ell_2)(\Delta(f))$; its filtration is defined
by $(\cO^\circ)_{\leq n} = \{\ell\in\cO^* | \ell(\m^{n+1}) = 0\}$.  

Note that we have a category equivalence ${\bf FSHA} \to {\bf LCA}$, 
where ${\bf LCA}$ is the category of Lie coalgebras, taking
$\cO$ to $\m/\m^2$, equipped with the cobracket induced by 
$\Delta - \Delta^{2,1}$. Then the composed functor 
${\bf LCA} \to {\bf FSHA} \to {\bf FilAlg}$ is $\c\mapsto U(\c^*)$
(recall that $\c^*$ is a Lie algebra). 

$(\cO_{\g^*},\Delta_0) = (\wh S(\g),\Delta_0)$ is a graded FSH algebra. 
Its restricted dual is the graded algebra $S(\g^*)$. Recall that 
$\cO_{\g^*}$ is also a Poisson algebra. We define the set of 
Poisson traces on $\cO_{\g^*}$ as the 
subspace of all $\ell\in\cO_{\g^*}^\circ$, such that $\ell(\{u,v\}) = 0$
for any $u,v\in\cO_{\g^*}$. Then $\{$Poisson traces on $\cO_{\g^*}\}
\subset \cO_{\g^*}^\circ$ identifies with $S(\g^*)^\g \subset S(\g^*)$; 
this is a graded subalgebra of $\cO_{\g^*}^\circ$. This defines a graded 
algebra structure on $\{$Poisson traces on $\cO_{\g^*}\}$. 

Consider the FSH algebra $(\cO_{\g^*},{}^\rho\Delta_0)$. It is 
isomorphic (as a filtered vector space) to  $(\cO_{\g^*},\Delta_0)$, 
and this isomorphism induces an algebra isomorphism between their associated
graded FSH algebras. It follows that we have an isomorphism 
of filtered vector spaces 
between the filtered algebra $(\cO_{\g^*},{}^\rho\Delta_0)^\circ$
and $S(\g^*)$, and the associated graded of this morphism is
an algebra isomorphism $\on{gr} ((\cO_{\g^*},{}^\rho\Delta_0)^\circ) \to 
S(\g^*)$. 

Recall that the vector spaces underlying $(\cO_{\g^*},\Delta_0)^\circ$ and 
$(\cO_{\g^*},{}^\rho\Delta_0)^\circ$ are the same (i.e., $\cO_{\g^*}^\circ$). 
We claim that the canonical inclusion $\{$Poisson traces on $\cO_{\g^*}\}
\subset (\cO_{\g^*},{}^\rho\Delta_0)^\circ$ is a morphism of filtered 
algebras. 
Indeed, let us denote by $\cdot_\rho$ (resp., $\cdot$) the product of 
$(\cO_{\g^*},{}^\rho\Delta_0)^\circ$ (resp., $(\cO_{\g^*},\Delta_0)^\circ$). 
Let $\ell_1,\ell_2$ be Poisson 
traces on $\cO_{\g^*}$. Then for any $x\in\cO_{\g^*}$, we have 
$(\ell_1 \cdot_\rho \ell_2)(x) = (\ell_1 \otimes \ell_2)
(\rho \star \Delta_0(f) \star (-\rho))$. Now Leibniz's rule implies that 
$(\ell_1\otimes \ell_2)(\{u,v\})=0$ for any $u,v\in\cO_{\g^*}^{\wh\otimes 2}$, 
therefore 
$(\ell_1 \cdot_\rho \ell_2)(x) = (\ell_1\otimes \ell_2)(\Delta_0(x))
= (\ell_1\cdot \ell_2)(x)$. So $\{$Poisson traces on $\cO_{\g^*}\}
\subset (\cO_{\g^*},{}^\rho\Delta_0)^\circ$ is an algebra morphism. 
Since the filtrations on the vector spaces underlying
$(\cO_{\g^*},\Delta_0)^\circ$ and $(\cO_{\g^*},{}^\rho\Delta_0)^\circ$
are the same, and since the filtration on $\{$Poisson traces on $\cO_{\g^*}\}$
is induced by that of $(\cO_{\g^*},\Delta_0)^\circ$, this morphism is filtered,
and its associated graded is the canonical inclusion 
$S(\g^*)^\g \subset S(\g^*)$. 

Now the FSH algebra isomorphism $\cO_{G^*} \simeq 
(\cO_{\g^*},{}^\rho\Delta_0)$ (Corollary \ref{coroprinc})
induces a filtered algebra isomorphism 
$(\cO_{\g^*},{}^\rho\Delta_0)^\circ \to \cO_{G^*}^\circ = U(\g^*)$.
The fact that the associated graded of this morphism is the canonical 
isomorphism $S(\g^*) \to \on{gr}(U(\g^*))$ follows from the fact that the 
completed graded of the FSH algebras $\cO_{G^*}$ and 
$(\cO_{\g^*},{}^\rho\Delta_0)$ are both $(\cO_{\g^*},\Delta_0)$. 

We now compose the filtered algebra morphism 
$\{$Poisson traces on $\cO_{\g^*}\}
\subset (\cO_{\g^*},{}^\rho\Delta_0)^\circ$
with the filtered algebra isomorphism 
$(\cO_{\g^*},{}^\rho\Delta_0)^\circ \to \cO_{G^*}^\circ$
and obtain a filtered algebra morphism $S(\g^*)^\g \to U(\g^*)$, 
whose associated graded is the canonical inclusion $S(\g^*)^\g 
\subset S(\g^*)$.  

The situation may be summarized as follows: 
$$
\begin{matrix}
 & & (\cO_{\g^*},\Delta_0)^\circ = S(\g^*) & & \\ 
 & \scriptstyle{(a)} \nearrow & & & & \\ 
S(\g^*)^\g = \{\hbox{Poisson\ traces\ on\ }\cO_{\g^*}\} & & 
\scriptstyle{(c)}\uparrow & & & \\ 
 & \scriptstyle{(b)}\searrow & & & \\ 
 & & (\cO_{\g^*},{}^\rho\Delta_0)^\circ & \stackrel{(d)}{\to} & 
 \cO_{G^*}^\circ = U(\g^*) 
\end{matrix}
$$
Here $S(\g^*)^\g$ and $S(\g^*)$ are graded algebras, 
$(\cO_{\g^*},{}^\rho\Delta_0)^\circ$ and $\cO_{G^*}^\circ$
are filtered algebras; $(a)$ is a morphism of graded algebras, $(c)$ is 
an isomorphism of filtered vector spaces, $(b)$ and $(d)$ are morphisms 
of filtered algebras ($(d)$ is an isomorphism). The associated graded 
of $(c)$ is an isomorphism of graded algebras.    
\hfill \qed \medskip 

\begin{remark} The restricted dual of the isomorphism $\cO_{\g^*} \to 
\cO_{G^*}$ appearing in the above proof is an 
isomorphism of filtered vector spaces $\sigma : S(\g^*) \to U(\g^*)$, whose 
associated graded is the canonical isomorphism 
$S(\g^*) \to \on{gr}(U(\g^*))$. These properties are also satisfied by 
the symmetrization map $\on{Sym}$, however $\sigma$ depends on $\rho$, 
so in general $\on{Sym}$ and $\sigma$ are different.  
\end{remark} 

\begin{remark} One can check that the morphism $\theta$ is independent 
on the choice of $(\rho,\varphi)$ (these choices are described in Remark
\ref{rem:sigma} and in Theorem \ref{lift2}). 
\end{remark}

\section{Duality of QUE and QFSH algebras} \label{sect:duality}

In this section, we recall some facts from \cite{Dr:QG} (proofs can be found in
\cite{Gav}). Let us denote by ${\bf QUE}$ the category of quantized universal
enveloping (QUE) algebras and by ${\bf QFSH}$ the  category of quantized formal
series Hopf (QFSH) algebras. We denote by ${\bf QUE}_{\on{fd}}$ and
${\bf QFSH}_{\on{fd}}$ the subcategories corresponding to finite dimensional 
Lie bialgebras. 

We have contravariant functors ${\bf QUE}_{\on{fd}} \to {\bf QFSH}_{\on{fd}}$, 
$U\mapsto U^*$ and ${\bf QFSH}_{\on{fd}} \to {\bf QUE}_{\on{fd}}$, 
$\cO\mapsto \cO^\circ$. These functors are inverse to each other. 
$U^*$ is the full topological dual of $U$, i.e., the space of all  
continuous (for the $\hbar$-adic topology) $\KM[[\hbar]]$-linear maps
$U \to \KM[[\hbar]]$. 
$\cO^\circ$ the space of continuous $\KM[[\hbar]]$-linear forms 
$\cO\to \KM[[\hbar]]$, 
where $\cO$ is equipped with the $\m$-adic topology (here $\m\subset \cO$
is the maximal ideal). 

We also have covariant functors ${\bf QUE} \to {\bf QFSH}$, $U\mapsto U'$
and ${\bf QFSH} \to {\bf QUE}$, $\cO\mapsto \cO^\vee$. There functors are
also inverse to each other. $U'$ is a subalgebra of $U$, while $\cO^\vee$
is the $\hbar$-adic completion of $\sum_{k\geq 0} \hbar^{-k} \m^k \subset 
\cO[1/\hbar]$. 

We also have canonical isomorphisms $(U')^\circ \simeq (U^*)^\vee$
and $(\cO^\vee)^* \simeq (\cO^\circ)'$. 

If $\a$ is a finite dimensional Lie bialgebra and $U = U_\hbar(\a)$
is a QUE algebra quantizing $\a$, then $U^* = \cO_{A,\hbar}$ is a 
QFSH algebra quantizing the Poisson-Lie group $A$ (with Lie bialgebra $\a$), 
and $U' = \cO_{A^*,\hbar}$ is a QFSH algebra quantizing the Poisson-Lie 
group $A^*$ (with Lie bialgebra $\a^*$). If now $\cO = \cO_{A,\hbar}$
is a QFSH algebra quantizing $A$, then $\cO^\circ = U_\hbar(\a)$ is a 
QUE algebra quantizing $\a$ and $\cO^\vee = U_\hbar(\a^*)$ is a QFSH algebra 
quantizing $\a^*$. 

We now compute these functors explicitly in the case of cocommutative  
QUE and commutative QFSH algebras. If $U = U(\a)[[\hbar]]$ with 
cocommutative coproduct
(where $\a$ is a Lie algebra), then $U'$ is a completion of
$U(\hbar \a[[\hbar]])$; this is a flat deformation of $\wh S(\a)$
equipped with its linear Lie-Poisson structure. If $G$ is a formal group 
with function ring $\cO_G$, then $\cO := \cO_G[[\hbar]]$ is a QFSH algebra, 
and $\cO^\vee$ is a commutative QUE algebra; it is a quantization of 
$(S(\g^*)$, commutative product, cocommutative coproduct, co-Poisson 
structure induced by the Lie bracket of $\g)$. 

\section{Relation between twist quantization and its functional version}
\label{sect:rel}

Let us define a twist quantization of 
the coboundary Lie bialgebra $(\g,r,Z)$ as 
a pair $(J,\Phi)$, $J\in U(\g)^{\otimes 2}[[\hbar]]$,
$\Phi \in U(\g)^{\otimes 3}[[\hbar]]$, such that $\Phi$ 
is invariant, and $(J,\Phi)$ satisfies the 
twisted cocycle relation
\begin{equation}
\label{eq1}
J^{1,2}J^{12,3}=J^{2,3}J^{1,23}\Phi,  
\end{equation}
$(\varepsilon\otimes\id)(J)=(\id\otimes\varepsilon)(J)=1$,
$J = 1^{\otimes 2} + O(\hbar)$, $\Phi = 1^{\otimes 3} + O(\hbar)$, 
$\on{Alt}((J-1^{\otimes 2})/\hbar) = r + O(\hbar)$, 
$\on{Alt}((\Phi-1^{\otimes 3})/\hbar^2) 
= Z + O(\hbar)$. 
These conditions imply that $\Phi$ satisfies the pentagon relation, as well as 
$\eps^{(i)}(\Phi) = 1^{\otimes 2}$, $i=1,2,3$. 
(We know that such a twist quantization always exists
when $\g$ is triangular or quasi-triangular.) 
Our purpose is to relate twist quantization with its functional version.  

The first step is to show that $(J,\Phi)$ can be transformed into an 
admissible pair, in a sense which we now precise.  

\begin{definition}
\label{admissible}
{\it An element $x$ in a QUE algebra  
$U$ is admissible if $x\in 1 + \hbar U$, and  
if $\hbar \log x$ is in $U' \subset U$.} 
\end{definition}

We will use the isomorphism $U(\g)^{\otimes k}[[\hbar]] \simeq 
U(\g^{\oplus k})[[\hbar]]$ to view  $U(\g)^{\otimes k}[[\hbar]]$
as a QUE algebra. 

\begin{proposition}
\label{propadmi}
Any twist quantization $(J,\Phi)$ of a coboundary Lie bialgebra $(\g,r,Z)$ 
is gauge equivalent to an admissible twist quantization $(J',\Phi')$
(i.e., such that $J'$ and $\Phi'$ are admissible).
\end{proposition}

{\em Proof.} Let us set $U = U(\g)[[\hbar]]$. 
According to Proposition \ref{prop:assoc}, one can find 
an invariant $F \in U^{\wh\otimes 2}$, such that 
$F \in 1^{\otimes 2} + \hbar U_0^{\wh\otimes 2}$ and  
$\Phi' := {}^F\Phi = F^{2,3} F^{1,23} \Phi (F^{1,2}F^{12,3})^{-1}$
is admissible. In particular, $\Phi' \in 1^{\otimes 3} 
+ \hbar^2 U_0^{\wh\otimes 3}$.  

Then if we set $J_0 := JF$, we have 
$J_0^{1,2} J_0^{12,3} = J_0^{2,3} J_0^{1,23} \Phi'$, 
and $J_0 \in 1^{\otimes 2} + \hbar U_0^{\wh\otimes 2}$. 
For any $u\in 1^{\otimes 3} + \hbar U_0$, 
${}^uJ_0 := u^1 u^2 J_0 (u^{12})^{-1}$
is such that $({}^uJ_0,\Phi')$ is a twist quantization of 
$(\g,r,Z)$. It remains to find $u$ such that $J' := {}^u J_0$
is admissible. 

We will construct $u$ as a product $\cdots u_2 u_1$, where 
$u_n\in 1 + \hbar^n U_0$, in such a way that if 
$J_n := {}^{u_n\cdots u_1}J_0$, then 
$\hbar\log(J_n) \in U_0^{\prime\wh\otimes 2} + \hbar^{n+2} U_0^{\wh\otimes 2}$.  

We have already $\hbar\log(J_0) \in \hbar^2 U_0^{\wh\otimes 2}$. 

Expand $J_0 = 1^{\otimes 2} + \hbar j_1 + \cdots$, then $\on{Alt}(j_1) = r$. 
Moreover, the coefficient of $\hbar$ in $J_0^{1,2}J_0^{12,3} = 
J_0^{2,3}J_0^{1,23}\Phi$ yields $d(j_1) = 0$, where $d : U(\g)_0^{\otimes 2}
\to U(\g)_0^{\otimes 3}$ is the co-Hochschild differential. 
It follows that for some $a_1\in U(\g)_0$, we have $j_1 = r+d(a_1)$. 
Then if we set $u_1 := \exp(\hbar a_1)$ and $J_1 = {}^{u_1}J_0$, 
we get $J_1 \in 1^{\otimes 2} + \hbar r + \hbar^2 U_0^{\wh\otimes 2}$.   
Then $\hbar\log(J_1) \in \hbar^2 r + \hbar^3 U_0^{\wh\otimes 2}
\subset U_0^{\prime\wh\otimes 2} + \hbar^3 U_0^{\wh\otimes 3}$. 

Assume that for $n\geq 2$, we have constructed $u_1,\ldots,u_{n-1}$ such that 
$\alpha_{n-1} := \hbar\log(J_{n-1}) \in U_0^{\prime\wh\otimes 2} 
+ \hbar^{n+1} U_0^{\wh\otimes 2}$. 

Let us denote by $\bar\alpha$ the image of the class of 
$\alpha_{n-1}$ in $U(\g)_0^{\otimes 2} / (U(\g)_0^{\otimes 2})_{\leq n+1}$
under the isomorphism of this space with $(U_0^{\prime\wh\otimes 2} 
+ \hbar^{n+1} U_0^{\wh\otimes 2}) / (U_0^{\prime\wh\otimes 2} 
+ \hbar^{n+2} U_0^{\wh\otimes 2})$ (see Lemma \ref{lemma:quot}). 
Let $\alpha\in U(\g)_0^{\otimes 2}$ be a representative of $\bar\alpha$, then 
$\alpha_{n-1} = \alpha' + \hbar^{n+1}\alpha$, where $\alpha'\in 
U_0^{\prime\wh\otimes 2} + \hbar^{n+2} U_0^{\wh\otimes 2}$. Let us set  
$\varphi' := \hbar\log(\Phi')$, then the twist equation gives  
$$
(-\alpha'-\hbar^{n+1}\alpha)^{1,23} \star_\hbar 
(-\alpha'-\hbar^{n+1}\alpha)^{2,3} \star_\hbar 
(\alpha'+\hbar^{n+1}\alpha)^{1,2} \star_\hbar
(\alpha'+\hbar^{n+1}\alpha)^{12,3} = \varphi', 
$$
where $\star_\hbar$ is defined as in Appendix \ref{app:A}.
According to Lemma \ref{lemma:approx}, the image of this equality 
in $(U^{\wh\otimes 3} + \hbar^{n+1} U^{\prime\wh\otimes 3}) / 
(U^{\wh\otimes 3} + \hbar^{n+2} U^{\prime\wh\otimes 3})  
\simeq U(\g)^{\otimes 3} / (U(\g)^{\otimes 3})_{\leq n+1}$
is $d(\bar\alpha)$, where $d$ is the co-Hochschild differential on  
$U(\g)_0^{\otimes \cdot} / (U(\g)_0^{\otimes\cdot})_{\leq n+1}$. 
Since $n\geq 2$, the relevant cohomology group vanishes, so 
$\bar\alpha = d(\bar\beta)$, where $\bar\beta\in U(\g)_0 
/(U(\g)_0)_{\leq n+1}$. Let $\beta\in U(\g)_0$ be a representative  
of $\bar\beta$ and set $u_n := \exp(\hbar^n\beta)$, $J_n := {}^{u_n}J_{n-1}$, 
$\alpha_n := \hbar\log(J_n)$. Then 
$$
\alpha_n = (\hbar^{n+1}\beta)^1 \star_\hbar (\hbar^{n+1}\beta)^2 
\star_\hbar \alpha_{n-1} \star_\hbar (-\hbar^{n+1}\beta)^{12}. 
$$
According to Lemma \ref{lemma:approx}, the image of $\alpha_n$
in 
$$
(U_0^{\wh\otimes 2} + \hbar^{n+1} U_0^{\prime\wh\otimes 2}) 
/ (U_0^{\wh\otimes 2} + \hbar^{n+2} U_0^{\prime\wh\otimes 2}) 
\simeq U(\g)_0^{\otimes 2}/(U(\g)^{\otimes 2}_0)_{\leq n+1}
$$ 
is $\bar\alpha - d(\bar\beta)=0$. So $\alpha_n$ belongs to 
$U_0^{\wh\otimes 2} + \hbar^{n+2} U_0^{\prime\wh\otimes 2}$, as 
required. This proves the induction step. 
\hfill \qed \medskip 

If now $(J',\Phi')$ is an admissible twist quantization, then 
$\rho := \hbar \log(J')_{|\hbar=0}$ and $\varphi := 
\hbar \log(\Phi')_{|\hbar=0}$
are formal functions on $\m_{\g^*}^{\wh\otimes 2}$ and 
$\m_{\g^*}^{\wh\otimes 3}$, solutions of the functional twist equation.

\section{Quantization of $\cO_G^\g \subset \cO_G$} \label{sect:O}

Using a (non necessarily admissible) twist quantization, we construct
a formal noncommutative deformation of the inclusion of algebras 
of Lemma \ref{lemma:O}:

\begin{proposition}
We have an injective algebra morphism 
$\cO_G^\g[[\hbar]] \hookrightarrow \cO_{G,\hbar}$
deforming $\cO_G^\g\subset \cO_G$, where 
$\cO_{G,\hbar}$ is a quantization of the PFSH algebra 
$\cO_G$ and $\cO_G^\g[[\hbar]]$ is the trivial deformation of 
the commutative algebra $\cO_G^\g$ (it is also commutative).  
\end{proposition}

{\em Proof.} Let us first construct the QFSH algebra $\cO_{G,\hbar}$. 
For $x\in U(\g)[[\hbar]]$, set 
${}^J\Delta_0(x) = J\Delta_0(x) J^{-1}$, 
where $\Delta_0$ is the usual cocommutative coproduct. 
Then 
$U_\hbar(\g) = (U(\g)[[\hbar]],m_0,{}^J\Delta_0)$ is a quantization 
of the Lie bialgebra $\g$ (here $m_0$ is the product on $U(\g)$). 
The dual $\cO_{G,\hbar} := U_\hbar(\g)^*$ of this QUE algebra 
is a QFSH algebra quantizing the PFSH algebra $\cO_G$. 
The product in this QFSH algebra is defined by 
$(f \star g)(x) = (f \otimes g)(J\Delta_0(x)J^{-1})$
for $f,g\in U_\hbar(\g)^*$ and $x\in U_\hbar(\g)$. 

On the other hand, the FSH algebra $\cO_G$ is equal to $U(\g)^*$, 
and its product is defined by $(fg)(x) = (f\otimes g)(\Delta_0(x))$
for $f,g\in U(\g)^*$ and $x\in U(\g)$. 

We say that $f\in U(\g)^*$ is a trace iff $f(xy) = f(yx)$
for any $x,y\in U(\g)$. Then the inclusion
$\{$traces on $U(\g)\} \subset U(\g)^*$ 
identifies with $\cO_G^\g \subset \cO_G$. 
In the same way, we define $\{$traces on $U(\g)[[\hbar]]\}$; 
this is a subalgebra of $U(\g)[[\hbar]]^* \simeq \cO_G[[\hbar]]$, 
which identifies with $\cO_G^\g[[\hbar]]$. 

The canonical map $\{$traces on $U(\g)[[\hbar]]\} \to U_\hbar(\g)^*$
is an algebra morphism. Indeed, if $f_1,f_2$ are traces on $U(\g)[[\hbar]]$, 
then $f_1\otimes f_2$ is a trace on $U(\g)^{\otimes 2}[[\hbar]]$, so 
$(f_1 \star f_2)(x) = (f_1\otimes f_2)(J\Delta_0(x) J^{-1}) = 
(f_1\otimes f_2)(\Delta_0(x)) = (f_1 f_2)(x)$ for any $x\in U(\g)[[\hbar]]$,
so $f_1\star f_2 = f_1 f_2$. So we have obtained an algebra morphism 
$\cO_G^\g[[\hbar]] \to U_\hbar(\g)^* = \cO_{G,\hbar}$. It is clearly a
deformation of the canonical inclusion $\cO_G^\g\subset \cO_G$. 
\hfill \qed \medskip 

\section{Quantization of $S(\g^*)^\g \hookrightarrow U(\g^*)$}
\label{sect:SG}

Assume now that $(J,\Phi)$ is an admissible twist quantization. 
We will construct a formal deformation of the inclusion of 
algebras of Theorem \ref{theoprinc}. 

\begin{theorem}
There is an injective algebra morphism:
$$
\theta_\hbar : S(\g^*)^\g[[\hbar]] \hookrightarrow U_\hbar(\g^*), 
$$
where $U_\hbar(\g^*)$ is a quantization of $\g^*$.
Its reduction modulo $\hbar$ coincides with the morphism  
$S(\g^*)^\g \hookrightarrow U(\g)$ from Theorem \ref{theoprinc}. 
\end{theorem}

{\em Proof.} Recall that $U(\g)[[\hbar]]'$ is a cocommutative QFSH 
algebra; we denote by $m_0$, $\Delta_0$ its product and 
coproduct. 

Since $(\eps\otimes\id)(J) = (\id\otimes\eps)(J)=1$, we have 
$\hbar\log(J) \in \m_0^{\wh\otimes 2}$, where $\m_0\subset 
U(\g)[[\hbar]]'$ is the kernel of the counit. According to \cite{EH1}, 
Proposition 3.1, this implies that the inner automorphism 
$z\mapsto JzJ^{-1}$ of $U(\g)^{\otimes 2}[[\hbar]]$ restricts to an 
automorphism of $U(\g)^{\otimes 2}[[\hbar]]'$. 

We then equip $U(\g)[[\hbar]]'$ with the coproduct 
${}^J\Delta : x\mapsto J \Delta_0(x) J^{-1}$. 
Then $(U(\g)[[\hbar]]',m_0,{}^J\Delta_0)$ is a QFSH 
algebra. Its classical limit is the PFSH algebra 
$(\cO_{\g^*},m_0,P,{}^\rho\Delta_0)$. We have seen that 
this PSFH algebra is isomorphic to $\cO_{G^*}$, hence 
$(U(\g)[[\hbar]]',m_0,{}^J\Delta_0)$ is a quantization of $\cO_{G^*}$. 

It now follows from Section \ref{sect:duality} that 
$(U(\g)[[\hbar]]',m_0,{}^J\Delta_0)^\circ$ is a 
quantization of $U(\g^*)$, which we denote by $U_\hbar(\g^*)$. 

Let us say that $\varphi\in (U(\g)[[\hbar]]')^\circ$ is a trace if 
$\varphi(xy) = \varphi(yx)$ for any $x,y\in U(\g)[[\hbar]]'$. 
Then $\{$traces on $U(\g)[[\hbar]]'\} \subset (U(\g)[[\hbar]]')^\circ$
is a subalgebra. Indeed, if $\ell_1,\ell_2$ are traces then 
$\ell_1\otimes \ell_2$ is also a trace, so for 
$x,y\in U(\g)[[\hbar]]'$, we have $(\ell_1\ell_2)(xy) = 
(\ell_1\otimes \ell_2)(\Delta(x)\Delta(y)) = 
(\ell_1\ell_2)(yx)$. This inclusion 
identifies with the inclusion $(\cO_G[[\hbar]]^\vee)^\g \subset
\cO_G[[\hbar]]^\vee$. Indeed, the Drinfeld 
functors have the property that $(U')^\circ = (U^*)^\vee$ for any QUE algebra
$U$. 

Now we show that the map 
$\{$traces on $U(\g)[[\hbar]]'\} \subset (U(\g)[[\hbar]]',{}^J\Delta_0)^\circ$
is also an algebra morphism. Indeed, let $\cdot_J$ be the product of the latter
algebra. If $\ell_1,\ell_2$ are traces
and $x,y\in U(\g)[[\hbar]]$, then 
$(\ell_1 \cdot_J \ell_2)(x) = (\ell_1\otimes \ell_2)(J\Delta_0(x)J^{-1})
= (\ell_1\otimes \ell_2)(\Delta_0(x)) = (\ell_1\ell_2)(x)$, so 
$\ell_1 \cdot_J \ell_2 = \ell_1 \ell_2$. So we have constructed an 
algebra morphism $(\cO_G[[\hbar]]^\vee)^\g \to U_\hbar(\g)$. It is clearly a 
deformation of the morphism constructed in Theorem \ref{theoprinc}.   

Recall that $\cO_G[[\hbar]]^\vee$ is the $\hbar$-adic completion of 
$\sum_{k\geq 0} \hbar^{-k} \m_G^k \subset 
\cO_G((\hbar))$.\footnote{$\cO_G[[\hbar]]^\vee$ may be also be viewed 
as the formal Rees algebra associated to the decreasing filtration 
$\cO_G \supset \m_G \supset \m_G^2 \cdots$.} 
Then $\cO_G[[\hbar]]^\vee$ is a topologically 
free $\KM[[\hbar]]$-commutative algebra; its specialization at $\hbar=0$ is 
$\cO_G[[\hbar]]^\vee / \hbar \cO_G[[\hbar]]^\vee \simeq S(\g^*)$. 

The action of $\g$ on $\cO_G$ induces an action of $\g$ on 
$\cO_G[[\hbar]]^\vee$. 
Then $(\cO_G[[\hbar]]^\vee)^\g$ is the $\hbar$-adic completion of 
$\sum_{k\geq 0} \hbar^{-k} (\m_G^k)^\g$. 
We have an inclusion of topologically free $\KM[[\hbar]]$-algebras
$(\cO_G[[\hbar]]^\vee)^\g \subset \cO_G[[\hbar]]^\vee$. 

Now the dual of the symmetrization map induces an algebra isomorphism 
$\wh S(\g^*) = \cO_\g\simeq \cO_G$ (dual to the exponential map $\g\to G$). 
This isomorphism induces a $\g$-equivariant isomorphism of 
$\cO_G[[\hbar]]^\vee$ with the $\hbar$-adic completion of
$\sum_{k\geq 0} \hbar^{-k} \m_\g^k \subset \cO_\g((\hbar))$. 
So we have an algebra isomorphism $\cO_G[[\hbar]]^\vee \simeq S(\g^*)[[\hbar]]$. 
It restricts to an isomorphism $(\cO_G[[\hbar]]^\vee)^\g \simeq S(\g^*)^\g[[\hbar]]$. 

Composing its inverse with the morphism $(\cO_G[[\hbar]]^\vee)^\g \to 
U_\hbar(\g^*)$, we get the announced morphism
$S(\g^*)^\g[[\hbar]] \to U_\hbar(\g^*)$. 
\hfill \qed
\medskip 

\section{The quasitriangular case}
\label{sect:Sem}

A quasitriangular Lie bialgebra (QTLBA) is a pair $(\g,r')$, where 
$\g$ is a Lie algebra and $r'\in \g^{\otimes 2}$ is such that 
$\on{CYB}(r')=0$ and $t:= r'+r^{\prime 2,1} \in S^2(\g)^\g$. Any QTLBA 
gives rise to a coboundary Lie bialgebra $(\g,r,Z)$, where 
$r=(r'-r^{\prime 2,1})/2$ and $Z = [t^{1,2},t^{2,3}]/4$. 
We call a QTLBA {\it nondegenerate} if $\g$ is finite dimensional and
$t$ is nondegenerate. 

Let $D : \g^* \to \g^*$ be the composition of 
the Lie cobracket $\delta : \g^* \to \wedge^2(\g^*)$ with the Lie bracket 
of $\g^*$. It is a derivation and a coderivation, and it induces a 
derivation of $U(\g^*)$, which we also denote by $D$ (or 
sometimes $D_{\g^*}$). 

\begin{proposition}
For any scalar $s$, $C_s := \on{Ker}(\delta - s (D\otimes \id)\circ \Delta_0)$
is a commutative subalgebra of $U(\g^*)$. 
\end{proposition}

{\em Proof.} The condition $\ell\in C_s$ means that 
for any $u,v\in \cO_{G^*}$, we have $\ell(\{u,v\} - s D^*(u)v) = 0$
(here $D^*$ is the derivation of $\cO_{G^*}$ dual to the coderivation 
$D$). 

Let $\ell_1,\ell_2$ belong to $C_s$. Then for any $u,v\in \cO_{G^*}$, 
\begin{align*}
(\ell_1\ell_2)(\{u,v\} - s D^*(u)v) & 
= (\ell_1\otimes \ell_2)(\{\Delta(u),\Delta(v)\}
- s \Delta(D^*(u)) \Delta(v)) 
\\ & = (\ell_1\otimes \ell_2)(\{u^{(1)},v^{(1)}\} \otimes u^{(2)}v^{(2)} 
+ u^{(1)}v^{(1)} \otimes \{u^{(2)},v^{(2)}\} 
\\ & - s D^*(u^{(1)})v^{(1)} \otimes u^{(2)}v^{(2)} -  u^{(1)}v^{(1)} 
\otimes sD^*(u^{(2)})v^{(2)}) = 0,
\end{align*} hence $\ell_1\ell_2\in C_s$. 
Here $\Delta$ is the coproduct of $\cO_{G^*}$. 

Moreover, we constructed in \cite{EGH} an element $\varrho\in
\m_{G^*}^{\wh\otimes 2}$, such that $\Delta'(u) = \varrho \star \Delta(u) \star
(-\varrho)$ for any $u\in\cO_{G^*}$; if $(U_\hbar(\g),\cR)$ is 
any quantization of $(\g,r')$, then $\hbar\log(\cR) \in 
\m_\hbar^{\wh\otimes 2}$, where $\m_\hbar \subset U_\hbar(\g)'$
is the augmentation ideal, and the reduction of $\hbar\log(\cR)$ 
mod $\hbar$ equals $\varrho$. Then it follows from $(S^2\otimes S^2)
(\cR) = \cR$ that $(S_\cO^2\otimes S_\cO^2)(\hbar\log\cR) = \hbar\log\cR$, 
where $S$ is the antipode of $U_\hbar(\g)$ and 
$S_\cO = S_{|U_\hbar(\g)'}$ is the antipode of 
$U_\hbar(\g)'\subset U_\hbar(\g)$; 
since the specialization for $\hbar=0$ of $\hbar^{-1}(S_\cO^2 -\id)$
is $D^*$, we get 
$(D^* \otimes \id + \id \otimes D^*)(\varrho) = 0$. 

Then if $\ell_1,\ell_2\in C_s$, then 
$(\ell_2\ell_1)(u) = (\ell_1\otimes \ell_2)(\Delta'(u))
= (\ell_1\otimes \ell_2)(\varrho \star \Delta(u) \star (-\varrho))
= (\ell_1\otimes \ell_2)(\Delta(u)) + \sum_{n\geq 1} (1/n!)
(\ell_1\otimes \ell_2)(\{\varrho,\{\varrho,\ldots,\{\varrho,\Delta(u)\}\})$. 
Now if $f\in \cO_{G^*}^{\wh\otimes 2}$, then $(\ell_1\otimes \ell_2)
(\{\varrho,f\}) = s(\ell_1\otimes \ell_2)((D^*\otimes \id + 
\id \otimes D^*)(\varrho)f) = 0$. 
It follows that $\ell_2\ell_1 = \ell_1 \ell_2$. 
\hfill \qed \medskip

\begin{remark}
If $A$ is a quasitriangular Hopf algebra with antipode $S$, 
set $C_{s,A} := \{\ell\in A^* | \forall a,b\in A, 
\ell(ab) = \ell(bS^{-2s}(a))\}$ for any $s\in\ZZ$. Then it follows from 
\cite{Dr:coco} that $C_{s,A}$ is a commutative algebra, and that we have
isomorphisms $C_s \simeq C_{s+2}$ for any $s\in\ZZ$. The isomorphism 
takes $\ell\in C_s$ to $\overline\ell\in C_{s+2}$ defined by 
$\overline\ell(x) = \ell(xu^{-1}S(u))$, where $u = m\circ 
(\id\otimes S)(R)$ ($m,R$ are the product and $R$-matrix of $A$). 
The definition of 
$C_{s,A}$ can be generalized to $s\in\KM$ when $A = (U_\hbar(\g),\cR)$
is a quasitriangular QUE Hopf algebra. Define $U_\hbar(\g^*)$
as $(U_\hbar(\g)')^\circ = (U_\hbar(\g)^*)^\vee \supset U_\hbar(\g)^*$. 
Then $C_{s,\hbar} := \{\ell\in (U_\hbar(\g)')^\circ | \forall
a,b\in U_\hbar(\g)', \ell(ab) = \ell(b (S^2)^{-s}(a)) \}$
is a commutative subalgebra of $U_\hbar(\g^*)$, and its reduction 
modulo $\hbar$ is contained in $C_s$. In this case, $u^{-1}S(u)$ does 
not necessarily belong to $U_\hbar(\g)'$, therefore $C_{s,\hbar}$ 
and $C_{s+2,\hbar}$ are not necessarily isomorphic. 
\end{remark}

\begin{remark}
If $(\g,r,Z)$ is a coboundary Lie bialgebra, then $r$ is $D$-invariant iff 
$(\mu\otimes \id)(Z)$ is symmetric (where $\mu$ is the Lie bracket of $\g$).
Otherwise, if we set $\varrho := \rho^{2,1} \star (-\rho)$, then  
$(D^*\otimes \id + \id \otimes D^*)(\varrho) \neq 0$, so unless $s=0$, 
one cannot prove that $C_s$ is commutative. 
\end{remark}

For each nondegenerate QTLBA $(\g,r')$, Semenov-Tian-Shansky defined 
an algebra morphism  $\Theta : Z(U(\g)) \to U(\g^*)$, where $Z(A)$
denotes the center of an algebra $A$ (\cite{STS1}). 
Let us recall the construction of $\Theta$.
There are unique Lie algebra morphisms $L,R : \g^* \to \g$, 
defined by $L(\ell) = (\ell\otimes \id)(r')$, $R(\ell) = 
-(\id\otimes \ell)(r')$
for any $\ell\in\g^*$. We denote by $\alpha : U(\g^*) \to U(\g)$
the composed map $U(\g^*) \stackrel{\Delta_0}{\to} U(\g^*)^{\otimes 2}
\stackrel{L \otimes (S_0\circ R)}{\to} 
U(\g)^{\otimes 2} \stackrel{m_0}{\to} U(\g)$. Here $m_0,\Delta_0$
are the standard product and coproduct maps, we still denote by $L,R$
the algebra morphisms induced by $L,R$, and $S_0$ denotes the antipode of 
$U(\g)$. The associated graded of the map $\alpha$ is the isomorphism 
$S(\g^*) \to S(\g)$ induced by $t$, hence $\alpha$ is an isomorphism. 
Then $\Theta : Z(U(\g)) \to U(\g^*)$ is defined as the restriction of
$\alpha^{-1}$ to $Z(U(\g))$; one can prove that it is an algebra morphism.
 
We will show, together with Proposition \ref{8:10}: 

\begin{proposition} \label{prop:sem} \label{8:4}
$\on{Im}(\Theta) = C_1 \subset U(\g^*)$.  The associated graded of 
$C_1$ (for the degree filtration of $U(\g^*)$) is $S(\g^*)^\g$. 
\end{proposition}

\begin{remark} Let $\theta$ be as in Theorem \ref{theoprinc}. 
The image of $\theta : S(\g^*)^\g \to U(\g^*)$ 
is $\{$Poisson traces on $\cO_{G^*}\}$, i.e., 
this is $\on{Ker}(\delta)$, where $\delta : 
U(\g^*) \to \wedge^2 U(\g^*)$ is the co-Poisson map of $U(\g^*)$. 
So the images of $\Theta$ and $\theta$ do not coincide.
\hfill \qed \medskip 
\end{remark}

Let us now construct a deformation $\Theta_\hbar$ of 
$\Theta$. The following lemma is proved in \cite{Dr:coco}. 

\begin{lemma} \label{8:9}
Let $(A,\Delta,R)$ be a quasitriangular Hopf algebra
with antipode $S$. 
Define a linear map $\alpha_A : A^* \to A$ by $\alpha_A(\ell) = 
(\ell\otimes \id)(R^{21}R)$. Then $\alpha_A$ induces an 
algebra morphism $C_{1,A} \to Z(A)$. 
\end{lemma}

\begin{lemma} \label{lemma:p7}
Assume moreover that $A$ is finite dimensional and $R^{2,1}R$
is nondegenerate. Then the map $C_1 \to Z(A)$ is a linear
isomorphism. Its inverse induces an algebra morphism $\Theta_A : 
Z(A) \to A^*$. 
\end{lemma}

{\em Proof.} We have to check that if $\ell\in A^*$ is such that 
$\alpha_A(\ell) \in Z(A)$, then $\ell$ is a trace. The condition 
$\alpha_A(\ell) \in Z(A)$ means that 
for any $a\in A$, we have $(\ell \otimes \id)([R^{2,1}R,1\otimes a])=0$. 
It follows that for any $a\in A$, we have 
$S^{-1}(a^{(4)})
(\ell\otimes \id)([R^{2,1}R, a^{(2)}S^{-1}(a^{(1)}) \otimes a^{(3)}]) =0$. 
Since $R^{2,1}R$ commutes with the image of $\Delta_A$,   
$(\ell\otimes \id)((a^{(2)} \otimes S^{-1}(a^{(4)}) a^{(3)})
[R^{2,1}R, S^{-1}(a^{(1)}) \otimes 1]) 
 =0$. Therefore 
$(\ell\otimes \id)((a^{(2)} \otimes 1) R^{2,1}R (S^{-1}(a^{(1)}) \otimes \id)) 
= \eps(a) (\ell\otimes \id)(R^{2,1} R)$. 

Since $R^{2,1}R$ is nondegenerate, this means that for any $b\in A$, 
we have $\ell(a^{(2)} b S^{-1}(a^{(1)})) = \eps(a)\ell(b)$. Replacing 
$a\otimes b$ by $a^{(1)} \otimes S(a^{(2)}) b$, we get 
$\ell(bS^{-1}(a)) = \ell(S(a)b)$, so that $\ell\in C_1$.  
\hfill \qed \medskip  

The QUE algebra version of these lemmas is 1), 2) of the following
proposition. 
Let $(\g,r')$ be a QTLBA and let $(U_\hbar(\g),\Delta,\cR)$ be a 
quantization of $(\g,r')$. 

\begin{proposition} \label{8:10}
1) The linear map $U_\hbar(\g)^* \to U_\hbar(\g)$, $\ell\mapsto 
(\ell\otimes \id)(\cR\cR^{2,1})$ extends to a map 
$\alpha_\hbar : U_\hbar(\g^*) \to U_\hbar(\g)$. 

2) If $(\g,r')$ is nondegenerate, then $\alpha_\hbar$ is a linear 
isomorphism, and it restricts to an algebra isomorphism $C_{1,\hbar} \to 
Z(U_\hbar(\g))$. 

3) Proposition \ref{8:4} is true. 
\end{proposition}

{\em Proof.} Let us prove 1). 
Define $L_\hbar, R'_\hbar : U_\hbar(\g)^* \to U_\hbar(\g)$
by $L_\hbar(\xi) = (\xi\otimes \id)(\cR)$, $R'_\hbar(\xi) = 
(\id\otimes \xi)(\cR)$. According to \cite{EH1}, $\hbar\log(\cR)
\subset (U_\hbar(\g)'_0)^{\wh\otimes 2} \subset U_\hbar(\g)'_0 
\wh\otimes \hbar U_\hbar(\g)_0$, so that $\log(\cR) \in 
U_\hbar(\g)'_0 \wh\otimes U_\hbar(\g)_0$. According to \cite{EGH}, appendix, 
the image of $\log(\cR)$
in $(\m_{G^*}/\m_{G^*}^2) \wh\otimes U(\g)_0$  (by reduction mod $\hbar$
followed by projection) is $r'$. It follows that $\cR\in U_\hbar(\g)' 
\wh\otimes U_\hbar(\g)$, therefore $L_\hbar$ extends to a map 
$U_\hbar(\g^*) \to U_\hbar(\g)$; this map is necessarily a QUE algebra 
morphism. The quasitriangularity identities imply 
that the image of $\cR$ in $\cO_{G^*} \wh\otimes 
U(\g)$ has the form $\on{exp}(\rho)$, where $\rho\in \m_{G^*}\otimes \g$
is a lift of $r$. It follows that the reduction mod $\hbar$ of $L_\hbar$
is the morphism induced by $\g^*\to\g$, $\ell \mapsto (\ell\otimes\id)(r)$. 
In the same way, $R'_\hbar$ extends to a (anti)morphism $U_\hbar(\g^*) \to 
U_\hbar(\g)$. 

Define $\alpha_\hbar : U_\hbar(\g^*) \to U_\hbar(\g)$, by 
$x\mapsto m \circ (L_\hbar\otimes R'_\hbar) \circ \Delta$. 
Then $\alpha_\hbar$ extends $\ell\mapsto (\ell\otimes \id)(\cR\cR^{2,1})$. 

Let us prove 2). The reduction mod $\hbar$ of $\alpha_\hbar$ is $\alpha$, 
which is a linear isomorphism; hence $\alpha_\hbar$ is a linear 
isomorphism. The second part is proved as Lemma \ref{8:9}. 

Let us prove Proposition \ref{8:4}. Assume that $U_\hbar(\g)$ is as in
\cite{EK}, hence $U_\hbar(\g) \simeq U(\g)[[\hbar]]$ as algebras. Then
$Z(U_\hbar(\g)) \simeq Z(U(\g))[[\hbar]]$. 2) implies that $\alpha$ induces an
isomorphism $(\on{mod\ }\hbar)(C_{1,\hbar}) \to Z(U(\g))$; here $(\text{mod\
}\hbar)$ is the reduction modulo $\hbar$. On the other hand, 
$(\text{mod\ }\hbar)(C_{1,\hbar})\subset C_1$, therefore $\Theta(Z(U(\g)))
\subset C_1$. 

The map $\delta - (D\otimes \id)\circ \Delta_0 : U(\g^*) \to 
U(\g^*)^{\otimes 2}$ is filtered, and its associated graded is the dual 
$\delta : S(\g^*) \to \wedge^2(S(\g^*))$ of the Poisson bracket of 
$S(\g)$. We have a surjective morphism $S(\g)_\g  = S(\g) / \{\g,S(\g)\}
\twoheadrightarrow S(\g)/\{S(\g),S(\g)\}$ to the cokernel of this Poisson
bracket, hence $\on{Ker}(\delta) \hookrightarrow (S(\g)_\g)^* = S(\g^*)^\g$. We
have $\on{gr}(C_1) \subset \on{Ker}(\delta)$, hence $\on{gr}(C_1) \subset 
S(\g^*)^\g$. Now since $\Theta$ is filtered and its associated graded takes 
$\on{gr}(Z(U(\g))) \simeq S(\g)^\g$ to $S(\g^*)^\g$, we get $\on{gr}(C_1)
= S(\g^*)^\g$ and $\Theta(Z(U(\g))) = C_1$. 
\hfill \qed \medskip

We denote by $\Theta_\hbar : 
Z(U_\hbar(\g)) \to U_\hbar(\g^*)$ the algebra morphism 
inverse to $\alpha_\hbar$. $\Theta_\hbar$
is the QUE algebra version of $\Theta_A$ defined above. 

The image of $\Theta_\hbar$ is $C_{1,\hbar}$.
When the quantization is an in \cite{EK}, $U_\hbar(\g) \simeq 
U(\g)[[\hbar]]$, so this image is not the same as that of $\theta_\hbar$, 
which is $\{$traces on $U_\hbar(\g)'\} = C_{0,\hbar}$. 
Therefore in this case,  
the images of $\theta_\hbar$ and $\Theta_\hbar$ do not coincide.

\section{On the canonical derivation of $\cO_{G^*}$} \label{sect:D}

Let $(\a,\mu_\a,\delta_\a)$ be a finite dimensional Lie bialgebra. 
Then $\cO_A$ is a Poisson-Lie group, dual to $U(\a)$. 
Set $D_\a := \mu_\a \circ \delta_\a$, then $D_\a$ is a 
derivation of $U(\a)$, such that if $U_\hbar(\a)$
is any quantization of $U(\a)$ with antipode $S$, then 
$D_\a = \hbar^{-1}(S^2-\id)_{|\hbar=0}$ (see \cite{Dr:coco}). 
It follows that the dual derivation $D_\a^*$ of $\cO_A$ has 
the same property. 

When $\a = (\g,r')$ is a quasitriangular Lie bialgebra,
$D_\a$ is inner, given by $D_\a(x) = -[\mu(r'),x]$ for any $x\in U(\g)$; 
here $\mu$ is the Lie bracket of $\g$ (see \cite{Dr:coco}). 

\begin{proposition} \label{prop:inner}
If $\g$ is a nondegenerate quasitriangular Lie bialgebra, then
the derivation $D_{\g^*}^*$ of $\cO_{G^*}$ is inner, i.e., 
there exists a function $h\in \cO_{G^*}$ such that 
$D_{\g^*}^*(f) = \{h,f\}$ for any $f\in \cO_{G^*}$.   
\end{proposition}

{\em Proof.} We assume that $\g$ is the double $\a_+\oplus \a_-$
of a Lie bialgebra $\a_+$ (here $\a_- = \a_+^*$); 
the general case is similar. Then $\g^*$ is (as a Lie algebra) 
the direct sum $\a_+\oplus \a_-$. Let $A_\pm$ be the formal groups 
corresponding to $\a_\pm$. The morphism $\alpha : U(\g^*) \to U(\g)$
is now $U(\a_+) \otimes U(\a_-) \to U(\g)$, $x_+ \otimes x_- \mapsto 
x_- S(x_+)$. The dual morphism $\alpha^* : \cO_{G} \to \cO_{G^*}$
takes $F\in \cO_G$ to $f\in \cO_{G^*}$ given by $f(g_+,g_-) := 
F(g_- g_+^{-1})$. 

\begin{lemma}
Let $D^*_\g$, $D^*_{\g^*}$ be the canonical derivations of 
$\cO_G$ and $\cO_{G^*}$. Then $\alpha^* \circ D^*_{\g} = 
D^*_{\g^*} \circ \alpha^*$. Moreover, $D^*_{\g} = 
\LL_{\mu(r')} - \RR_{\mu(r')}$, where $\mu$ is the Lie bracket of 
$\g$ and $\LL_a f(g) = (d/d\eps)_{|\eps=0} F(e^{\eps a}g)$, 
$\RR_a f(g) = (d/d\eps)_{|\eps=0} F(ge^{\eps a})$.  
\end{lemma}

{\em Proof of Lemma.} $D_{\g^*}$ is a coderivation, so $\Delta_0 : 
U(\g^*) \to U(\g^*)^{\otimes 2}$ intertwines $D_{\g^*}$ and 
$D_{\g^*} \otimes \id + \id\otimes D_{\g^*}$; $L$ and $R$ are 
Lie bialgebra morphisms, so they intertwine $D_{\g^*}$ and $D_{\g}$; 
$S$ commutes with $D_{\g}$; and $D_\g$ is a derivation, so $m_0$ intertwines 
$D_{\g} \otimes \id + \id \otimes D_\g$ with $D_\g$. 
Hence $\alpha \circ D_{\g^*} = D_\g \circ \alpha$. The first part follows. 

According to \cite{Dr:coco}, $D_\g(x) = -[\mu(r'),x]$, which implies the 
second part. \hfill \qed \medskip 

In \cite{STS2}, the image of the Poisson bracket on $G^*$ under the 
formal isomorphism $\alpha : G^*\to G$ dual to $\alpha^*$ was 
computed. Let $f,h\in \cO_{G^*}$ and 
$F = (\alpha^*)^{-1}(f)$, $H = (\alpha^*)^{-1}(h)$, then 
\begin{align} \label{PB:sem}
(\alpha^*)^{-1}(\{f,h\})(g) & = 
\langle (d_{\RR} - d_\LL) F(g) \otimes d_\LL H(g), r' \rangle
+ \langle (d_{\RR} - d_\LL) F(g) \otimes d_\RR H(g), (r')^{2,1} \rangle
\\ & \nonumber 
= \langle (d_\LL - d_\RR) F(g), L(d_\RR H(g)) - R(d_\LL H(g)) 
\rangle 
\end{align}
where $g\in G$, $d_\LL F(g), d_\RR F(g) \in \g^*$ are the left and right
differentials defined by 
$\langle d_\LL F(g) , a\rangle = (\LL_a F)(g)$, 
$\langle d_\RR F(g) , a\rangle = (\RR_a F)(g)$ 
for any $a\in \g$. 

\begin{lemma} \label{9:3}
There exists a function $H(g)\in \cO_G$ such that 
$L(d_{\RR} H(g)) - R(d_\LL H(g)) = \mu(r')$. 
\end{lemma}

{\em Proof of lemma.} We prove this when $\g$ is the double 
$\a_+ \oplus \a_-$ of a Lie bialgebra $\a_+$. Then set $a = (a_+,a_-)$
where $a_\pm\in \a_\pm$. We have $\g^* = \a_+ \oplus \a_-$, and we 
should solve: $d_{\RR} H_a(g) = \mu(r')_- + u_+(g)$, 
$d_{\LL} H_a(g) = \mu(r')_+ + u_-(g)$, 
where $u_\pm(g)$ are functions $G \to \a_\pm$.
Now $d_\LL H(g) = \Ad(g) (d_\RR H(g))$, hence 
$\mu(r')_+ + u_-(g) = \Ad(g)(\mu(r')_- + u_+(g))$. 
Let us decompose 
$g = g_- g_+^{-1}$, where $g_\pm \in A_\pm = \exp(\a_\pm)$, we get 
$\Ad(g_+^{-1})(u_+(g)) - \Ad(g_-^{-1})(u_-(g)) = \Ad(g_-^{-1})( \mu(r')_+) 
- \Ad(g_+^{-1})(\mu(r')_-)$. 
Therefore 
$$
u_+(g) = \Ad(g_+) \big( \Ad(g_-^{-1})(\mu(r')_+) 
- \Ad(g_+^{-1})(\mu(r')_-)\big)_+
$$
and the condition is
$$
d_\RR H(g)  = \Ad(g_+) \big( \Ad(g_+^{-1})(\mu(r')_-)\big)_-   
+ \Ad(g_+) \big( \Ad(g_-^{-1})(\mu(r')_+) \big)_+ , 
$$
i.e., 
\begin{equation} \label{flatness}
\RR_\alpha H_a(g) = 
\langle \mu(r')_-, \Ad(g_+)\big((\Ad(g_+^{-1})(\alpha))_+\big) \rangle 
+ \langle \mu(r')_+, \Ad(g_-) \big( (\Ad(g_+^{-1})(\alpha))_- \big) \rangle
\end{equation} 
for any $\alpha\in \g$. 
Let us denote by $A_\alpha(g)$ the r.h.s. of (\ref{flatness}). 

Let us compute $\RR_\alpha A_\beta - \RR_\beta A_\alpha$, 
for $\alpha,\beta\in\g$. Recall that $g = g_- g_+^{-1}$, then  
we have $\RR_\alpha(g) = g \alpha$, so $\RR_\alpha(g_\pm^{-1}) = 
\pm (\Ad(g_\pm^{-1})(\alpha))_\pm g_\pm^{-1}$. After computations, we find: 
$$
\RR_\alpha A_\beta - \RR_\beta A_\alpha = A_{[\beta,\alpha]} 
+ B_{\alpha,\beta}, 
$$ 
where
\begin{align*}
& B_{\alpha,\beta}(g) = - \langle [(\Ad^*(g_+^{-1})(\beta))_+,
(\Ad^*(g_+^{-1})(\alpha))_+], 
(\Ad(g_+^{-1})(\mu(r')_-))_-\rangle 
\\ & + \langle [(\Ad^*(g_+^{-1})(\beta))_-,(\Ad^*(g_+^{-1})(\alpha))_-], 
(\Ad(g_-^{-1})(\mu(r')_+))_+\rangle . 
\end{align*}
Now for $u,v\in\a_+$, we have 
\begin{align*}
& \langle [u,v], (\Ad(g_+^{-1})(\mu(r')_-))_- \rangle 
= \langle [u,v], \Ad(g_+^{-1})(\mu(r')_-) \rangle
\\ & 
= \langle [\Ad^*(g_+)(u), \Ad^*(g_+)(v)], \mu(r')_- \rangle
= \langle [\Ad^*(g_+)(u), \Ad^*(g_+)(v)], \mu(r') \rangle
\\ & 
= \langle \Ad^*(g_+)(u) \otimes \Ad^*(g_+)(v), 
\delta(\mu(r')) \rangle =0, 
\end{align*}
since $\delta(\mu(r')) = 0$ (see \cite{Dr:coco}). In the same way, the second
term of $B_{\alpha,\beta}(g)$ vanishes. Hence the system 
(\ref{flatness}) has a solution
(it is unique if we impose that $H$ vanishes at the origin). 
\hfill \qed \medskip 

{\em End of proof of Proposition \ref{prop:inner}.} 
Now if $h = -\alpha^*(H)$
with $H$ as in Lemma \ref{9:3} and for any $f\in\cO_{G^*}$, we have 
\begin{align*}
& D_{\g^*}(f) = \alpha^*(D_\g^*(F)) = 
\alpha^*((\LL_{\mu(r')} - \RR_{\mu(r')})(F))
\\ & = \alpha^*(\langle (d_\LL - d_\RR)(F)(g), R(d_\LL H(g)) - L(d_\RR H(g))
\rangle) = \{h,f\}. 
\end{align*}
\hfill \qed \medskip

\appendix

\section{Proof that associators can be made admissible} 
\label{app:A}

In \cite{EH2}, Proposition 3.2, 2) should read ``Assume that $x\in U'$ and
for any trees $R$, etc.". This affects Proposition 4.5 in \cite{EH2}, 
because the proof implicitly relies on the statement of 
Proposition 3.2 of \cite{EH2} without the assumption $x\in U'$. 
Below we prove a particular case of Proposition 4.5 of \cite{EH2} 
(the general case is similar). 

\begin{proposition} \label{prop:assoc}
Let $\g$ be a Lie algebra and let 
$\Phi \in U(\g)^{\otimes 3}[[\hbar]]$ be an invariant solution of 
the pentagon equation, such that $\eps^{(i)}(\Phi)=1^{\otimes 2}$ 
for $i=1,2,3$, $\Phi = 1^{\otimes 3} + O(\hbar)$ and 
$\on{Alt}(\Phi) = O(\hbar^2)$. Then there exists 
an invariant twist $F\in U(\g)^{\otimes 2}[[\hbar]]$, such that 
$\eps^{(i)}(F)=1$, $i=1,2$, $F=1^{\otimes 2} + O(\hbar)$, and 
${}^F\Phi = F^{2,3}F^{1,23} \Phi (F^{1,2}F^{12,3})^{-1}$ is admissible, 
i.e., $\hbar\log({}^F\Phi) \in (U(\g)[[\hbar]]')^{\wh\otimes 3}$.  
\end{proposition}  

{\em Proof.} We will construct $F$ as a product $\cdots F_2F_1$, 
where $F_n$ belongs to 
$1^{\otimes 2} + \hbar^i U_0^{\wh\otimes 2}$ and is such that if 
$\Phi_n := {}^{F_n\cdots F_1}\Phi$, then 
$\hbar\log(\Phi_n) \in (U'_0)^{\wh\otimes 2} + \hbar^{n+2}
U_0^{\wh\otimes 3}$. Here $U = U(\g)[[\hbar]]$ and the index $0$ denotes 
the augmentation ideals. 

We first construct $F_1$. Expand $\Phi = 1^{\otimes 3} + \hbar \phi_1 + \cdots$, 
then $d(\phi_1) = 0$ and $\on{Alt}(\phi_1)=0$, hence $\phi_1 = d(\psi_1)$, 
where $\psi_1\in (U_0^{\wh\otimes 2})^\g$ (here $d$ is the co-Hochschild 
differential of $(U_0^{\otimes \cdot})^\g$). 
We then set $F_1 = 1^{\otimes 2} + \hbar \psi_1$; 
we get $\Phi_1=1^{\otimes 3}+\hbar^2\phi'_2+\cdots$. Then $\hbar\log(\Phi_1) 
\in \hbar^3 U_0^{\wh\otimes 3}$. 

Now $d(\phi'_2)=0$, so there exists $\psi_2\in (U_0^{\wh\otimes 2})^\g$
such that $\phi'_2 = Z + d(\psi_2)$, where $Z\in \wedge^3(\g)[[\hbar]]$. 
Set $F_2 := 1^{\otimes 2}+\hbar^2\psi_2$, we get $\hbar\log(\Phi_2) 
\in \hbar^3 Z + \hbar^4 U_0^{\wh\otimes 3} \subset (U'_0)^{\wh\otimes 3}
+ \hbar^4 U_0^{\wh\otimes 3}$. 

Let $n\geq 3$. Assume that we have constructed $F_1,\ldots,F_{n-1}$ 
and let us construct $F_n$. By assumption, 
$\Phi_{n-1}\in 1^{\otimes 3} + \hbar U_0^{\wh\otimes 3}$ is such that 
$\varphi_{n-1} := \hbar\log(\Phi_{n-1}) \in (U'_0)^{\wh\otimes 3} 
+ \hbar^{n+1} U_0^{\wh\otimes 3}$. 

\begin{lemma} \label{lemma:quot}
The quotient $(U'+\hbar^n U) / (U'+\hbar^{n+1}U)$
identifies with $U(\g)/U(\g)_{\leq n}$. In the same way, the quotient 
$(U_0^{\prime\wh\otimes k} + \hbar^n U_0^{\wh\otimes k}) 
/ (U_0^{\prime\wh\otimes k} + \hbar^{n+1} U_0^{\wh\otimes k})$ identifies with
$U(\g)_0^{\otimes k} / (U(\g)_0^{\otimes k})_{\leq n}$ 
and the quotient of $\g$-invariant subspaces 
$(U_0^{\prime\wh\otimes k} + \hbar^n U_0^{\wh\otimes k})^\g 
/ (U_0^{\prime\wh\otimes k} + \hbar^{n+1} U_0^{\wh\otimes k})^\g$
identifies with
$(U(\g)_0^{\otimes k})^\g / (U(\g)_0^{\otimes k})^\g_{\leq n}$. 
\end{lemma}

The inverse of the first isomorphism takes the class of $\beta\in U(\g)$
to the class of $\hbar^n\beta\in U'+\hbar^n U$. 
Let $\bar\alpha\in (U(\g)_0^{\otimes 3})^\g/
(U(\g)_0^{\otimes 3})_{\leq n+1})^\g$ be the image of the class of 
$\varphi_{n-1}$ under the above isomorphism. Let $\alpha\in 
(U(\g)_0^{\otimes 3})^\g$ be a representative of $\bar\alpha$. Then we have 
$\varphi_{n-1} = \varphi + \hbar^{n+1}\alpha$, where 
$\varphi\in (U'_0)^{\wh\otimes 3} + \hbar^{n+2} U_0^{\wh\otimes 3}$. 
Now $\varphi_{n-1}$ satisfies the pentagon equation, so we get 
\begin{align} \label{star:id}
(-\varphi-\hbar^{n+1}\alpha)^{1,2,34} \star_\hbar 
(\varphi+\hbar^{n+1}\alpha)^{2,3,4} 
\star_\hbar (\varphi+\hbar^{n+1}\alpha)^{1,23,4} 
& \star_\hbar (\varphi+\hbar^{n+1}\alpha)^{1,2,3} 
\\ & \nonumber 
\star_\hbar (-\varphi-\hbar^{n+1}\alpha)^{12,3,4} = 0,  
\end{align}
where $a\star_\hbar b$ is the CBH product for the Lie bracket $[a,b]_\hbar = 
[a,b]/\hbar$. 

\begin{lemma} \label{lemma:approx}
Assume that $n\geq 2$. 
If $f_1,f_2\in (U'_0)^2 + \hbar^{n+1}U_0$ and $g,h\in \hbar^n U_0$, 
then $(f_1 + g) \star_\hbar (f_2 + h) = g + h$ modulo 
$(U'_0)^2 + \hbar^{n+1}U_0$. 
\end{lemma}

{\em Proof of Lemma.} The contribution of the degree $1$ 
part of the CBH series is $f_1+g + f_2 + h$, which gives $g+h$
modulo $(U'_0)^2 + \hbar^{n+1}U_0$. 

We now prove that $[(U'_0)^2 + \hbar^n U_0, (U'_0)^2 + \hbar^n U_0]_\hbar 
\subset (U'_0)^2 + \hbar^{n+1} U_0$. Indeed, we have 
$[U'_0,U'_0]_\hbar \subset U'_0$, hence
$[(U'_0)^2,(U'_0)^2]_\hbar \subset (U'_0)^2$;  
 $[\hbar^n U_0,\hbar^n U_0]_\hbar \subset \hbar^{2n-1}U_0 
\subset \hbar^{n+1} U_0$; and 
$[U'_0,\hbar^n U_0]_\hbar \subset \hbar^n U_0$ 
since $U'_0 \subset \hbar U_0$, so that 
$[(U'_0)^2,\hbar^n U_0]_\hbar \subset \hbar^n (U_0 U'_0 + U'_0 U_0) 
\subset \hbar^{n+1} (U_0)^2$ again because $U'_0 \subset \hbar U_0$.  

It follows that the contributions of all the higher degree parts of the CBH
series belong to $(U'_0)^2 + \hbar^{n+1}U_0$. This implies the lemma. 
\hfill \qed \medskip

{\em End of proof of Proposition \ref{prop:assoc}.} 
Lemma \ref{lemma:approx} implies that the image of (\ref{star:id}) in 
$(U^{\prime\wh\otimes 4} + \hbar^{n+1}U^{\wh\otimes 4}) / 
(U^{\prime\wh\otimes 4} + \hbar^{n+2}U^{\wh\otimes 4})
= U(\g)^{\otimes 4} / (U(\g)^{\otimes 4})_{\leq n+2}$
gives $d(\bar\alpha) = 0$, where 
$$
d : U(\g)^{\otimes 3} / (U(\g)^{\otimes 3})_{\leq n+2}
\to U(\g)^{\otimes 4} / (U(\g)^{\otimes 4})_{\leq n+2}
$$
is the map induced by the co-Hochschild differential. 

According to \cite{Dr:QH}, the cohomology of the complex 
$C^2\to C^3 \to C^4$ vanishes, where 
$C^k = (U(\g)_0^{\otimes k})^\g/(U(\g)_0^{\otimes k})_{\leq n+2}^\g$. 

It follows that there exists $\bar\beta\in C^2$, such that 
$\bar\alpha = d(\bar\beta)$. Let $\beta\in (U(\g)_0^{\otimes 2})^\g$
be a representative of $\bar\beta$. Set $F_n := \exp(\hbar^n\beta)$
and $\Phi_n = {}^{F_n}\Phi_{n-1}$. 

We get 
$$
\varphi_n = f_n^{2,3}\star_\hbar f_n^{1,23}
\star_\hbar \varphi_{n-1} \star_\hbar (-f_n^{12,3})\star_\hbar (-f_n^{1,2}), 
$$
where $f_n = \hbar^{n+1}\beta$. According to Lemma \ref{lemma:approx}, 
the class of $\varphi_n$ in 
$((U'_0)^{\wh\otimes 3} + \hbar^{n+1} U_0^{\wh\otimes 3}) / 
((U'_0)^{\wh\otimes 3} + \hbar^{n+2} U_0^{\wh\otimes 3})$ is 
$\bar\alpha - d(\bar\beta) = 0$, hence 
$\varphi_n \in (U'_0)^{\wh\otimes 3} + \hbar^{n+2} U_0^{\wh\otimes 3}$. 
This proves the induction step. 
\hfill \qed 

\end{document}